\documentclass{ifacconf}
\usepackage{amsmath}
\usepackage{float}
\usepackage{graphicx}      
\usepackage{natbib}        

\begin{document}
\setcitestyle{authoryear}
\pdfminorversion=4
\begin{frontmatter}

\title{Modelling the Incomplete Intermodal Terminal Location Problem } 

\author[First]{Oudani, M.} 

\address[First]{International University of Rabat, FIL, TIC Lab \\
Rocade Rabat-Sale, Morocco (e-mail: mustapha.oudani@uir.ac.ma).}

\begin{abstract}                
In this paper, we introduce and study the incomplete version of the intermodal terminal location problem. It's a generalization of the classical version by relaxing the assumption that the induced graph by located terminals is complete. We propose a mixed integer program to model the problem and we provide several extensions. All models are tested through validation in CPLEX solver. Numerical results are reported using well-known data set from the literature. 
\end{abstract}


\end{frontmatter}

\section{Introduction}
Due to its reliability, sustainability and its economical competitiveness, Intermodal Transportation (IT) has gained a good reputation. In fact, despite its lack of flexibility in the transport chain, intermodal transportation operators strikes to respect time schedules. Furthermore, the intermodal transportation is gaining ground over road transportation due to large investments in infrastructure development. Moreover, several customers  focus, today, on environmentally solutions in the transport industry. For all these reasons, intermodal transportation has attracted the attention of researchers and industry operators. The location of intermodal terminals is among the most challenging issues in the scientific literature. We study in this paper the Intermodal Terminal Location Problem in incomplete networks. The remainder of this article is organized as follows: we provide the state of the art in the section 2, problem description is given in the section 3, we propose several extensions of the original formulation in the section 4 and we conclude in the section 5.

\section{State of the art}
The scientific literature on Intermodal Location Problems is relatively recent, but the number of articles dealing with this subject is steadily increasing. \cite{arnold01} modeled the intermodal transshipment centers location as a linear program and proposed several extensions of the basic version. \cite{artmann03} studied the sustainability issues of the traffic shift from road to rail.  \cite{bontekoning04} provided a review synthesis about intermodal transportation in the field of operations research. \cite{limbourg09} proposed an iterative heuristic based on the p-median problem and on the multimodal assignement problem. \cite{sorensen12} proved that the intermodal terminal location problem is NP-hard and proposed efficient heuristics to solve it. \cite{tsamboulas07} developed a methodology for the policy measures assessment for modal shift to intermodal transportation. \cite{lin14} proposed a simplified version of the Sorensen model and proposed two efficient math-heuristics to solve it. \cite{oudani14} solved the problem using a genetic algorithm and proposed a new intermodal cost evaluation.  \cite{lin16} proposed a two-stage matheuristic approach to solve the problem. They proposed a program reformulation of the problem and test it using randomly generated data set. \cite{abbassi17} proposed a bi-objective model for transportation of agricultural products from Morocco to
Europe and developed heuristics to solve it. Recently, \cite{mostert18} proposed a bi-objective mathematical model minimizing the transportation and environmental costs objectives. \cite{abbassi18}  studied the robust  intermodal freight transport problem and proposed two solutions approaches for solving the problem. To the best of our knowledge, the current paper is the first to consider the terminals network incompleteness.  

\section{Problem description}

\subsection{The incomplete version motivation}
As in hub location problems, most studies in intermodal terminal location problems assume a complete inter-terminals network, that is, every terminal pair is interconnected. In the incomplete network studied here, we relax this assumption. In fact, assuming a complete network increase the total investment cost and connecting all terminals directly may also become unnecessary and expensive. 

\begin{figure}
\begin{center}
\includegraphics[width=8.4cm]{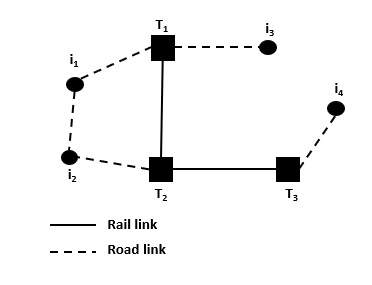}   
\caption{Small incomplete network} 
\label{fig:incom}
\end{center}
\end{figure}

 Fig.~\ref{fig:incom} shows a small incomplete intermodal terminal network with 3 terminals and 4 customers. The induced graph by terminals is not complete. For instance, there is no rail link between terminal $T_1$ and $T_3$. 

\subsection{Mathematical formulation}

Let be the following parameters:\\
$I:$ set of customers. \\
$K:$ set of potentials sites for intermodal terminals. \\
$l: $ number of links between located terminals. \\
$q_{ij}: $ the amount of goods to be transported from the customer $i$ to customer $j$. \\
$c_{ij}^{km}:$ the unit cost for intermodal transportation from customer $i$ to customer $j$ through the two terminal $k$ and $m$. \\
$c_{ij} :$ the unit unimodal cost for routing goods from customer $i$ to customer $j$. \\
$f_k:$ the cost for location of the terminal $k$. \\
$C_k:$ the capacity of the terminal $k$. \\
Let be the following decision variables: \\
$z_{km} :$ if the inter-terminals link between $k$ and $m$ is established, 0 otherwise. \\
$w_{ij} :$ the amount of goods transported by road from the customer $i$ to customer $j$ \\
$x_{ij}^{km}$ the amount of goods transported by road from the customer $i$ to customer $j$ through the two terminal $k$ and $m$ \\
The incomplete intermodal terminal location problem may be modeled as follows: 
\begin{equation}\label{eq:eq1}
\begin{split}
    \text{ Min } \sum_{i, j \in I }\sum_{k, m \in K} c_{ij}^{km} x_{ij}^{km} + \\
    \sum_{i, j \in I}c_{i j}w_{i j}+ \sum_{k \in K} f_k z_{kk}
    \end{split}
\end{equation}
 Subject to : 
\begin{equation}\label{eq:eq2}
    \sum_{k, m \in K} x_{i j}^{km}+ w_{ij}=q_{ij}, \forall i, j \in I
\end{equation}
\begin{equation}\label{eq:eq3}
    \sum_{i,j \in I} \sum_{ m \in K}x_{i j}^{km}+ \sum_{i,j \in I} \sum_{ m \in K}x_{i j}^{mk} \leq C_k z_{kk}, \forall k \in K
\end{equation}
\begin{equation}\label{eq:eq4}
    z_{km} \leq z_{kk}, \forall k, m \in K
\end{equation}
\begin{equation}\label{eq:eq5}
    z_{km} \leq z_{mm}, \forall k, m \in K
\end{equation}
\begin{equation}\label{eq:eq6}
    z_{km} = z_{mk} , \forall k, m \in K
\end{equation}
\begin{equation}\label{eq:eq7}
    \sum_{k,m \in K} z_{km}=l
\end{equation}
\begin{equation}\label{eq:eq8}
   x_{i j}^{km} \leq z_{km} , \forall k,m \in K, \forall i,j \in I
\end{equation}
The objective function~(\ref{eq:eq1}) minimizes the total cost for routing goods by road and by using the intermodal rail-road transportation. The constraint~(\ref{eq:eq2}) states that the sum of the amount routed directly from $i$ to $j$ and that through the terminals $k$ and $m$ is equal exactly the total goods to be transported form $i$ to $j$. The constraint~(\ref{eq:eq3}) guarantees the respect of the terminals capacities. Inequalities ~(\ref{eq:eq4}) and ~(\ref{eq:eq5}) state that an inter-terminal link is established if the two terminals are located. The equation ~(\ref{eq:eq6}) specify that if an inter-terminal link is established in the two directions. The equation ~(\ref{eq:eq7}) states the establishment of the given number of the inter-terminal links and the last equation ~(\ref{eq:eq8})  forbids the inter-terminal flows between closed terminals. The decision variable  $ z_{km}$  controls the inter-terminal links to open. Thus, if for some reason, a railway link is inconceivable between two terminals du to geographic, economic or environmental constraints, this link is prohibited by  $z_{kl}=0$. \\
The model is a Mixed 0-1 Integer Program (MIP). If we denote $p$ the cardinal of $K$ and by $n$ the cardinal of $I$, then the program has $n^2p^2+3p^2+n^2+p+1$ constraints and $n^2p^2+n^2+p^2$ variables.

\begin{prop}   
If the triangular inequality holds for units cost then $x_{ij}^{kk}=0, \forall i,j \in I, \forall k \in K$.
\end{prop}
This proposition demonstrates that this constraint used in several mathematical formulations in the literature is an unnecessary constraint [10].
\begin{pf}    
The intermodal unit cost is calcultaed as follows: $c_{i j}^{km}=c_{ik}+ \alpha c_{km} + c_{mj}$ where $c_{ik}, c_{km}, c_{mj}$ are respectively the cost between the customer $i$ and the terminal $k$, the inter-terminal cost between $k$ and $m$ and the cost between the terminal $m$ and customer $j$ with $\alpha <1$. The coefficient $\alpha$ is the discount parameter expressing the scale economy generated by using rail mode between the terminal $k$ and $m$. For instance, this coefficient is assumed to be equal to 0.5 in the work of S\" orensen et al. [10].  Since, $c_{i j}^{kk}=c_{ik}+  + c_{kj} \geq c_{ij}$, then $x_{ij}^{kk}=0$
\end{pf}
\begin{prop}
Let be $ |K|=p$ the cardinal of the potentials sites set. If $l > \frac{p(p-1)}{2}$ then the problem is unfeasible. 
\end{prop}
\begin{pf}
A complete graph (fully connected)  with $p$ vertices has at most $\frac{p(p-1)}{2}$ edges.
\end{pf}

\subsection{Exact solutions}
The model is validated by implementation in CPLEX 12.6 solver. To report numerical results, we used the instances randomly generated by S\"o rensen et al. Customers and potentials sites coordinates are randomly generated between $(0,0)$ and $(10^4,10^4)$. The goods demands are generated from the interval $[0,500]$. The investment cost $f_k$ in the interval $[0,5.10^5]$ and potentials sites capacities are drawn from $[0,10^4]$. After that, cost $c_{ij} $ is equal the euclidean distance between the customers $i$ and $j$ while $c_{i j}^{km}=c_{ik}+ \frac{1}{2} c_{km} + c_{mj}$. The exact solutions for some instances are reported in the  table~\ref{tb:exact1}. Optimal solutions for larger instances (more than 90 customers and 40 potential location) are not found in 1 hour. 

\begin{table}[ht]
\begin{center}
\caption{Exact solutions}\label{tb:exact1}
\begin{tabular}{cccc}
Instance & Cost $ (\times 10^7 )$ & Time (s) & \# terminals   \\\hline
\\
10C10L2TL  &   10,2       &   0,55   &  6 \\
10C10L4TL  &   9,62     &     1,09  &   6  \\
10C10L6TL  &   9,42    &      1,00   &  7  \\
10C10L8TL  &   9,36     &     0,69   &  8  \\
10C10L10TL  &  9,31   &       1,42  &   6  \\
10C10L12TL  &  9,29    &      1,06  &   9  \\

20C10L2TL   &  52,8     &     2,95  &   3  \\
20C10L4TL  &   52,7    &      2,52   &  5  \\
20C10L6TL  &   52,7  &        1,92  &   5  \\
20C10L8TL  &   52,7   &       1,63 &    5  \\
20C10L10TL  &  52,7    &      1,64   &  5  \\
20C10L12TL  &  52,8     &     3,08  &   6  \\

40C10L2TL   &  201,5    &     31,02 &   7  \\
40C10L4TL   &  201      &     22,55  &  8  \\
40C10L6TL   &  200,8    &     36,92 &   8  \\
40C10L8TL   &  200,7    &     18,45 &   8  \\
40C10L10TL  &  200,6    &     18,78 &   8  \\
40C10L12TL  &  200,6     &    13,83  &  8  \\

80C10L2TL  &   778,1    &     338,77 &   3 \\
80C10L4TL   &  778,8     &    1349,20 &  4  \\
80C10L6TL   &  778,8    &     599,92  &  4 \\
80C10L8TL & 779,8  &  1529,13 & 5 \\
80C10L10TL  & 779,8 &  891,28 &   5 \\
 \\\hline
\end{tabular}
\end{center}
\end{table}

\section{Extensions} 
\subsection{Minimizing the number of inter-terminals links }
Instead of minimizing the number of terminals to be located, we consider the problem of minimizing the number of links between a given $q$ number of terminals to be located. This problem can be modeled as follows: 
\begin{equation}\label{eq:eq9}
    \text{ Min }  \sum_{i, j \in I }\sum_{k, m \in K} c_{ij}^{km} x_{ij}^{km} + \sum_{i, j \in I}c_{ij}w_{i j}+\sum_{k, m \in K}c_{km}z_{km}
\end{equation}
Subject to: 
$$
    ~(\ref{eq:eq2})-~(\ref{eq:eq3})-~(\ref{eq:eq4})-~(\ref{eq:eq5})-~(\ref{eq:eq6})-~(\ref{eq:eq8})
$$
\begin{equation}\label{eq:eq10}
    \sum_{k \in K}z_{kk}=q
\end{equation}
We report in the table~\ref{tb:exact2} optimal solutions of this version for some instances with different values of number of terminals. 

\begin{table}[ht]
\begin{center}
\caption{Exact solutions}\label{tb:exact2}
\begin{tabular}{cccc}
Instance & Cost $ (\times 10^7)$ & Time (s) & \# links  \\\hline
\\
10C10L2T  &  11,24  &   2,31 &  2   \\
10C10L4T  &   10,2       &   0,55   &  5 \\
10C10L6T & 9,82  &  0,63  &    7\\
10C10L8T  &  9,74  &  0,45 &  9   \\
10C10L10T   &  9,74 &  0,44 &    10 \\
20C20L4T  &   50,9  &     161,41 &   6 \\
20C20L8T   &  48,2  &  270,94   &   12  \\
20C20L12T &  46,9  &   3442,39   &   17   \\
20C20L16T &  *  &   *   &   *   \\
20C20L20T &  *  &   *   &   *   \\
30C30L4T   & 117,2   &  1658,91  &  10**  \\
30C30L8T   & *   &  *  &  *  \\
30C30L12T   & *   &  *  &  *  \\
30C30L16T   & *   &  *  &  *  \\
30C30L20T   & *   &  *  &  *  \\
40C40L4T   & *   &  *  &  *  \\
40C40L8T   & *   &  *  &  *  \\
40C40L12T   & *   &  *  &  *  \\
40C40L16T   & *   &  *  &  *  \\
40C40L20T   & *   &  *  &  *  \\
50C50L8T   & *   &  *  &  *  \\
50C50L12T   & *   &  *  &  *  \\
50C50L16T   & *   &  *  &  *  \\
50C50L20T   & *   &  *  &  *  \\
60C60L4T   & *   &  *  &  *  \\
60C60L8T   & *   &  *  &  *  \\
60C60L12T   & *   &  *  &  *  \\
60C60L16T   & *   &  *  &  *  \\
60C60L20T   & *   &  *  &  *  \\
 \\\hline
\end{tabular}
\end{center}
\end{table}
\subsection{Minimizing the handling cost in terminals}
This version aims to minimize the operational handling cost in terminals. We denote $t_{km}$ the handling cost in terminal $k$ to terminal $m$ (we consider this cost as asymetric i.e $t_{km} \neq t_{mk}$). This problem may be formulated as follows: 
\begin{equation}\label{eq:eq11}
\begin{split}
    \text{ Min }  \sum_{i, j \in I }\sum_{k, m \in K}
    c_{ij}^{km} x_{ij}^{km} + \sum_{i, j \in I}c_{ij}w_{i j}+ \\
    \sum_{k, m \in K}(t_{km}+t_{mk})z_{km}+ \sum_{k \in K} f_k z_{kk}
\end{split}
\end{equation}
Subject to: 
$$
    ~(\ref{eq:eq2})-~(\ref{eq:eq3})-~(\ref{eq:eq4})-~(\ref{eq:eq5})-~(\ref{eq:eq6})-~(\ref{eq:eq7})-~(\ref{eq:eq8})
$$
We report in the table~\ref{tb:exact3} optimal solutions of this version for some instances with different values of number of inter-terminals links. 
\begin{table}[ht]
\begin{center}
\caption{Exact solutions}\label{tb:exact3}
\begin{tabular}{cccc}
Instance & Cost $ (\times 10^7 )$ & Time (s) & \#  links  \\\hline
\\
10C10L2TL  &  11,24  &   2,31 &  2   \\
10C10L4TL  &   10,2       &   0,55   &  5 \\
10C10L6TL & 9,82  &  0,63  &    7\\
10C10L8TL  &  9,74  &  0,45 &  9   \\
10C10L10TL   &  9,74 &  0,44 &    10 \\
20C20L4TL  &   50,9  &     161,41 &   6 \\
20C20L8TL   &  48,2  &  270,94   &   12  \\
20C20L12TL &  46,9  &   3442,39   &   17   \\
40C40L4TL   & *   &  *  &  *  \\
40C40L8TL   & *   &  *  &  *  \\
40C40L12TL   & *   &  *  &  *  \\
40C40L16TL   & *   &  *  &  *  \\
40C40L20TL   & *   &  *  &  *  \\
50C50L8TL   & *   &  *  &  *  \\
50C50L12T   & *   &  *  &  *  \\
50C50L16TL   & *   &  *  &  *  \\
50C50L20TL   & *   &  *  &  *  \\
60C60L4TL   & *   &  *  &  *  \\
60C60L8TL   & *   &  *  &  *  \\
60C60L12TL   & *   &  *  &  *  \\
60C60L16TL   & *   &  *  &  *  \\
60C60L20TL   & *   &  *  &  *  \\
 \\\hline
\end{tabular}
\end{center}
\end{table}
\subsection{Intermodal $(p,l)-$terminal location problem }
In this version both the number $p$ of the terminals to be located and the number $l$ of the inter-terminals links are given. This problem can be modeled as follows: 
\begin{equation}\label{eq:eq12}
    \text{ Min }  \sum_{i, j \in I }\sum_{k, m \in K} c_{ij}^{km} x_{ij}^{km} + \sum_{i, j \in I}c_{ij}w_{i j}
\end{equation}
Subject to: 
$$
    ~(\ref{eq:eq2})-~(\ref{eq:eq3})-~(\ref{eq:eq4})-~(\ref{eq:eq5})-~(\ref{eq:eq6})-~(\ref{eq:eq7})-~(\ref{eq:eq8})-~(\ref{eq:eq10})
$$
We report in the table~\ref{tb:exact4} optimal solutions of this version for some instances.
\begin{table}[ht]
\begin{center}
\caption{Exact solutions}\label{tb:exact4}
\begin{tabular}{ccc}
Instance & Cost $ (\times 10^7 )$ & Time (s)   \\\hline
\\
10C10L2T2TL  &  6,24  &   2,31   \\
10C10L4T4TL  &   9,12       &   0,35   \\
10C10L6T6TL & 9,82  &  0,53  \\
10C10L8T8TL  &  9,74  &  0,55   \\
10C10L10T10TL   &  9,74 &  0,47  \\
20C20L4T4TL  &   54,9  &     171,67  \\
20C20L8T8TL   &  45,2  &  250,65     \\
20C20L12T12TL &  41,9  &   3456,65      \\
40C40L4T4TL   & *   &  *   \\
40C40L8T8TL   & *   &  *    \\
40C40L12T12TL   & *   &  *   \\
40C40L16T16TL   & *   &  *   \\
40C40L20T20TL   & *   &  *    \\
50C50L8T8TL   & *   &  *   \\
50C50L12T12TL   & *   &  *   \\
50C50L16T16TL   & *   &  *   \\
50C50L20T20TL   & *   &  *    \\
60C60L4T4TL   & *   &  *   \\
60C60L8T8TL   & *   &  *    \\
60C60L12T12TL   & *   &  *   \\
60C60L16T16TL   & *   &  *   \\
60C60L20T20TL   & *   &  *    \\
 \\\hline
\end{tabular}
\end{center}
\end{table}

\section{Conclusion}

In this paper, we proposed a general version of the classical Intermodal Terminal Location Problem (ITLP) when the induced graph by located terminals is not necessarily complete. We formulate the problem as 0-1 linear program and proposed several extentions. We reported numerical results on data set instances given in the literature using CPLEX solver. As perspectives, we envision: 
\begin{enumerate}
    \item To develop efficient heuristics to solve larger instances
    \item To combine the problem with routing problems
    \item To study a hybrid hub intermodal terminal location problem considering incomplete networks 
\end{enumerate}



\end{document}